\documentclass[a4paper,11pt]{amsart}
\usepackage{amsmath,amssymb,amsthm}
\usepackage{enumitem}
\setlist[enumerate]{font=\normalfont,leftmargin=*,beginpenalty=10000}
\usepackage{xcolor}
\definecolor{refred}{RGB}{128,0,0}
\definecolor{citeblue}{RGB}{0,0,128}
\usepackage[colorlinks=true,citecolor=citeblue,linkcolor=refred,urlcolor=citeblue,pdftitle={Cube-freeness and a square formula for canonical trace ideals of Stanley--Reisner rings},pdfauthor={Sora Miyashita}]{hyperref}
\theoremstyle{plain}
\newtheorem{theorem}{Theorem}[section]

\newtheorem{lemma}[theorem]{Lemma}
\newtheorem{corollary}[theorem]{Corollary}

\newtheorem{maintheorema}{Theorem}
\newtheorem{maintheoremb}[maintheorema]{Theorem}
\newtheorem{introcorollary}{Corollary}

\theoremstyle{definition}
\newtheorem{example}[theorem]{Example}

\newtheorem{introquestion}{Question}
\theoremstyle{remark}
\newtheorem{remark}[theorem]{Remark}
\DeclareMathOperator{\tr}{tr}
\DeclareMathOperator{\Hom}{Hom}
\DeclareMathOperator{\Ext}{Ext}
\DeclareMathOperator{\core}{core}

\newcommand{\kk}{\Bbbk}
\newcommand{\NN}{\mathbb N}
\newcommand{\ZZ}{\mathbb Z}
\newcommand{\ee}{\mathbf e}
\newcommand{\one}{\mathbf1}
\newcommand{\mm}{\mathfrak m}
\author[S. Miyashita]{Sora Miyashita}
\address{Department of Pure and Applied Mathematics, The University of Osaka, Suita, Osaka 565-0871, Japan}
\email{u804642k@ecs.osaka-u.ac.jp}
\date{}
\title[Canonical trace ideals of Stanley--Reisner rings]{Cube-freeness and a square formula for canonical trace ideals of Stanley--Reisner rings}
\begin{document}
\begin{abstract}
Let $I$ be a monomial ideal in $S=\kk[x_1,\ldots,x_n]$, and set $R=S/I$. We prove that if $I$ is $\mathbf{a}$-bounded and $M$ is a positively $\mathbf{b}$-determined $R$-module, then $\tr_R(M)$ is $(\mathbf{a}+\mathbf{b})$-bounded. For the canonical module, this gives the bound $2\mathbf{a}$. In particular, the canonical trace ideal of every Stanley--Reisner ring is a cube-free monomial ideal. We also prove that, for a Stanley--Reisner ring satisfying $(S_2)$, the canonical trace ideal is the square of its radical if and only if the ring satisfies $(G_2)$.
\end{abstract}
\maketitle

\section{Introduction}\label{sec:introduction}

Let $S=\kk[x_1,\ldots,x_n]$ be a polynomial ring over a field $\kk$, and let $R=\kk[\Delta]=S/I_\Delta$ denote the \emph{Stanley--Reisner ring} of a simplicial complex $\Delta$ on $[n]=\{1,\ldots,n\}$. We write $\mm=(x_1,\ldots,x_n)R$ for its graded maximal ideal. For a finitely generated $R$-module $M$, the \emph{trace ideal} is
\[
\tr_R(M)=\sum_{\varphi\in\Hom_R(M,R)}\varphi(M).
\] We are concerned with the trace ideal of the canonical module $\omega_R$, which we call the \emph{canonical trace ideal} of $R$. If $R$ is Cohen--Macaulay, this ideal defines the non-Gorenstein locus (see \cite[Lemma~2.1]{HHS}).

We say that $R$ is \emph{locally Gorenstein} if $R_{\mathfrak p}$ is Gorenstein for every prime $\mathfrak p\ne\mm$. By \cite[Theorem~A]{MV}, a Cohen--Macaulay Stanley--Reisner ring is locally Gorenstein if and only if its canonical trace ideal is one of $R$, $\mm$, and $\mm^2$. If $\Delta$ is connected and $R$ is locally Gorenstein, then its canonical trace ideal is one of these three ideals even without the Cohen--Macaulay assumption (see \cite[Theorem~4.19]{KM26}).

We call a monomial ideal \emph{cube-free} if every exponent in its minimal monomial generators is at most two. We use the same definition for monomial ideals in a monomial quotient, considering only nonzero monomial generators. For a Stanley--Reisner ring $R=S/I_\Delta$, an ideal is cube-free if and only if its inverse image in $S$ is cube-free; we therefore omit the ambient polynomial ring when no confusion is possible. In particular, $R$, $\mm$, and $\mm^2$ are all cube-free. Moreover, if $\Delta$ is connected and $\dim R\geq3$, then $\tr_R(\omega_R)\ne\mm$ (see \cite[Corollary~4.20]{KM26}). Thus, in this range, the canonical trace ideal of a locally Gorenstein ring is either $R$ or $\mm^2$. In particular, $\tr_R(\omega_R)=\bigl(\sqrt{\tr_R(\omega_R)}\bigr)^2$. These observations lead to the following questions.

\begin{introquestion}\label{qu:cube-free}
Is the canonical trace ideal of every Stanley--Reisner ring cube-free?
\end{introquestion}

\begin{introquestion}\label{qu:square}
Under what hypotheses does $\tr_R(\omega_R)=\bigl(\sqrt{\tr_R(\omega_R)}\bigr)^2$ hold?
\end{introquestion}

Our first result gives an exponent bound for trace ideals over rings defined by monomial ideals. To state the result, give $S$ the multigrading $\deg x_i=\ee_i$, where $\ee_1,\ldots,\ee_n$ are the standard basis vectors of $\ZZ^n$. Set $\mathbf{x}=(x_1,\ldots,x_n)$, $\NN=\{0,1,\ldots\}$, and $\one=(1,\ldots,1)$. Inequalities between vectors are coordinatewise, and $[0,\mathbf{c}]$ denotes the set of $\mathbf{u}\in\NN^n$ with $\mathbf{u}\leq\mathbf{c}$.

For $\mathbf{c}\in\NN^n$, a monomial ideal in $S$ is \emph{$\mathbf{c}$-bounded} if its minimal generators divide $\mathbf{x}^{\mathbf{c}}$ (see \cite[Section~1]{AHZ22}). For a monomial ideal $I\subseteq S$, we use the same terminology for an ideal of $S/I$ whose minimal nonzero monomial generators have exponent vectors in $[0,\mathbf{c}]$. Thus the \emph{squarefree ideals} are precisely the $\one$-bounded ideals, and the cube-free ideals are precisely the $(2,\ldots,2)$-bounded ideals.

Let $\mathbf{c}\geq\one$. Following Miller \cite[Definition~2.1]{Mil00}, we call a finitely generated $\NN^n$-graded $S$-module $N$ \emph{positively $\mathbf{c}$-determined} if, for each $i$, multiplication by $x_i$ induces an isomorphism $N_{\mathbf{u}}\to N_{\mathbf{u}+\ee_i}$ whenever $u_i\geq c_i$. For $\mathbf{c}=\one$, this is the definition of a \emph{squarefree module} due to Yanagawa \cite[Definition~2.1]{Yan00}. When $R=S/I$, we call an $R$-module positively $\mathbf{c}$-determined or squarefree if it has the corresponding property as an $S$-module.

\begin{maintheorema}\label{thm:main-a}
Let $\mathbf{a},\mathbf{b}\in\NN^n$ with $\mathbf{a},\mathbf{b}\geq\one$. Let $I\subsetneq S$ be an $\mathbf{a}$-bounded monomial ideal, and set $R=S/I$. If $M$ is a positively $\mathbf{b}$-determined $R$-module, then $\tr_R(M)$ is $(\mathbf{a}+\mathbf{b})$-bounded. The canonical trace ideal $\tr_R(\omega_R)$ is $2\mathbf{a}$-bounded.
\end{maintheorema}

Theorem~\ref{thm:cube-free} proves the first assertion for any multigraded module generated in a translate of $[0,\mathbf{b}]$, and Corollary~\ref{cor:monomial-canonical-bound} gives the bound for $\omega_R$. Taking $\mathbf{a}=\one$ yields the following answer to Question~\ref{qu:cube-free}.

\begin{introcorollary}\label{cor:intro-cube-free}
The canonical trace ideal of every Stanley--Reisner ring is cube-free.
\end{introcorollary}

We next consider Question~\ref{qu:square}. Recall that $R$ \emph{satisfies $(S_2)$} if $\operatorname{depth}R_{\mathfrak p}\geq\min\{2,\dim R_{\mathfrak p}\}$ for every prime $\mathfrak p$. We say that $R$ \emph{satisfies $(G_r)$} if $R_{\mathfrak p}$ is Gorenstein whenever $\dim R_{\mathfrak p}\leq r$.

\begin{maintheoremb}\label{thm:main-b}
Let $R=\kk[\Delta]$ satisfy $(S_2)$. Then $R$ satisfies $(G_2)$ if and only if
\[
\tr_R(\omega_R)=\bigl(\sqrt{\tr_R(\omega_R)}\bigr)^2.
\]
\end{maintheoremb}

The radical in Theorem~\ref{thm:main-b} need not be maximal; see Example~\ref{ex:square-nonisolated}. If $\Delta=\core(\Delta)$ and $\dim\Delta\geq1$, then, under $(S_2)$, the equality holds precisely when $\Delta$ is a normal pseudomanifold (Corollary~\ref{cor:core-square}). Here the \emph{core} $\core(\Delta)$ is the subcomplex induced by the vertices that do not belong to every facet of $\Delta$ (see \cite[pp.~64--65]{Sta96}).

The proof of Theorem~\ref{thm:main-b} uses Gr\"abe's description of the canonical module to construct a pairing for normal pseudomanifolds (Lemma~\ref{lem:canonical-pairing}). Example~\ref{ex:square-hypotheses} shows that $(S_2)$ cannot simply be omitted and that $(G_2)$ cannot be replaced by $(G_1)$.

\section{Preliminaries}\label{sec:preliminaries}

In this section we fix notation and recall the facts about canonical modules and positively determined modules needed in the proofs. We then describe $(S_2)$ and $(G_2)$ in terms of links and recall Gr\"abe's description of the canonical module of a Stanley--Reisner ring. For background on monomial ideals and Stanley--Reisner rings, we refer to \cite{BH98,HH11}.

We retain the notation of the introduction. All modules are finitely generated. Unless otherwise stated, gradings and homogeneity refer to the $\ZZ^n$-grading. We use the shift convention $[M(\mathbf{a})]_{\mathbf{u}}=M_{\mathbf{u}+\mathbf{a}}$. The standard grading will also be used, with $\deg x_i=1$; in this grading, $\mm=\bigoplus_{j>0}R_j$.

For $F\subseteq[n]$, set $x_F=\prod_{i\in F}x_i$ and $\chi_F=\sum_{i\in F}\ee_i$. For $\mathbf{u}\in\NN^n$, the \emph{support} of $\mathbf{u}$, or of $\mathbf{x}^{\mathbf{u}}$, is $\operatorname{supp}(\mathbf{u})=\{i:u_i>0\}$.

For canonical modules of finite modules in the local and standard graded settings, see Schenzel \cite[Sections~1.2 and~3.3]{Sch98}. In our multigraded setting, let $R=S/I$ and let $M$ be a nonzero multigraded $R$-module of dimension $e$. Set $\omega_S=S(-\one)$. The \emph{canonical module} of $M$ is the $R$-module $\omega_M=\Ext_S^{n-e}(M,\omega_S)$, with the natural $\ZZ^n$-grading obtained from a multigraded free resolution of $M$. By multigraded local duality (see \cite[Corollary~6.1]{Mil00}), this is the graded $\kk$-dual of $H^e_{(x_1,\ldots,x_n)}(M)$. Here the \emph{graded dual} of a multigraded $\kk$-vector space $V$ is defined by $[V^\vee]_{\mathbf{u}}=\Hom_\kk(V_{-\mathbf{u}},\kk)$. We set $\omega_0=0$.

In the definition of a trace ideal, homomorphisms of all degrees are allowed. Thus a shift of a module does not change its trace ideal. Every homomorphism between multigraded modules is a finite sum of homogeneous homomorphisms. If $R=S/I$ with $I$ monomial, each nonzero multigraded component of $R$ is spanned by a monomial. It follows that the trace ideal of a multigraded $R$-module is a monomial ideal. We shall also use the fact that trace ideals commute with flat base change (see \cite[Lemma~1.5\,(iii)]{HHS}). A graded ring $R$ is \emph{quasi-Gorenstein} if $\omega_R\cong R(a)$ for some integer $a$, where the shift refers to the standard grading. For a Cohen--Macaulay ring with a canonical module, $\tr_R(\omega_R)=R$ if and only if $R$ is Gorenstein (see \cite[Lemma~2.1]{HHS}).

\begin{remark}\label{rem:positively-determined}
Let $\mathbf{a}\geq\one$. For a monomial ideal $I\subseteq S$, the following conditions are equivalent (see \cite[Example~2.4]{Mil00}):
\begin{enumerate}[label=\textup{(\roman*)}]
\item $I$ is $\mathbf{a}$-bounded;
\item $I$ is positively $\mathbf{a}$-determined;
\item $S/I$ is positively $\mathbf{a}$-determined.
\end{enumerate}
By \cite[Corollary~4.6]{Mil00}, a positively $\mathbf{a}$-determined module $N$ has a minimal multigraded free resolution with summands $S(-\mathbf{u})$, where $\mathbf{u}\in[0,\mathbf{a}]$. In particular, $N$ is generated in these multidegrees. Applying $\Hom_S(-,S(-\mathbf{a}))$ gives summands $S(-(\mathbf{a}-\mathbf{u}))$, which are again positively $\mathbf{a}$-determined. This property is preserved under kernels and cokernels of homogeneous maps of degree zero (see \cite[Remark~2.3]{Mil00}), so $\Ext_S^q(N,S(-\mathbf{a}))$ is positively $\mathbf{a}$-determined for every $q\geq0$.

If $I$ is $\mathbf{a}$-bounded and $\mathbf{c}\geq\mathbf{a}$, then a monomial ideal of $S/I$ is $\mathbf{c}$-bounded if and only if its inverse image in $S$ is $\mathbf{c}$-bounded. Indeed, the inverse image is generated by lifts of the minimal monomial generators of the ideal in $S/I$, together with the generators of $I$.
\end{remark}

A \emph{simplicial complex} on a finite set $V$ is a collection $\Delta$ of subsets of $V$ closed under taking subsets and containing $\{v\}$ for every $v\in V$. We call $V(\Delta)=V$ its \emph{vertex set}. Links and subcomplexes are regarded as complexes on their own vertex sets. The members of $\Delta$ are called \emph{faces}, and its maximal faces are called \emph{facets}. The \emph{dimension} of a face $F$ is $|F|-1$, and $\dim\Delta=\max_{F\in\Delta}(|F|-1)$ for $\Delta\ne\varnothing$. A complex is \emph{pure} if all its facets have the same dimension. In a pure $(d-1)$-dimensional complex, the $(d-2)$-faces are called \emph{ridges}.

For a complex on $[n]$, the \emph{Stanley--Reisner ideal} is $I_\Delta=(x_F:F\subseteq[n],\ F\notin\Delta)$. If $\Delta\ne\varnothing$, then $\dim\kk[\Delta]=\dim\Delta+1$. We allow both the empty collection $\Delta=\varnothing$ and the complex $\Delta=\{\varnothing\}$ on the empty vertex set $[0]=\varnothing$; their Stanley--Reisner rings are zero and $\kk$, respectively. We call $\Delta$ \emph{Cohen--Macaulay} over $\kk$ if $\kk[\Delta]$ is Cohen--Macaulay.

The \emph{link} of a face $F\in\Delta$ is $\operatorname{lk}_\Delta F=\{G\in\Delta:G\cap F=\varnothing,\ G\cup F\in\Delta\}$. If $\Delta_1$ and $\Delta_2$ have disjoint vertex sets, their \emph{join} is $\Delta_1*\Delta_2=\{F\cup G:F\in\Delta_1,\ G\in\Delta_2\}$. The \emph{suspension} of $\Delta$ is its join with a complex consisting of two isolated vertices.

A \emph{cone vertex} is a vertex that belongs to every facet. Let $C$ denote the set of cone vertices of $\Delta$. Following Stanley \cite[pp.~64--65]{Sta96}, we call $\core(\Delta)=\{F\in\Delta:F\cap C=\varnothing\}$ the \emph{core} of $\Delta$. Thus $\Delta=\core(\Delta)$ precisely when there are no cone vertices. Removing the cone vertices corresponds to removing polynomial variables: $\kk[\Delta]\cong\kk[\core(\Delta)][x_i:i\in C]$.

For $\Delta\ne\varnothing$, $\kk[\Delta]$ satisfies $(S_2)$ if and only if $\Delta$ is pure and $\operatorname{lk}_\Delta F$ is connected for every $F\in\Delta$ with $\dim\operatorname{lk}_\Delta F\geq1$ (see \cite[Section~1]{MT09}). Here and below, $F=\varnothing$ is allowed.

A pure complex is \emph{strongly connected} if any two facets $P$ and $Q$ can be joined by a sequence of facets $P=P_0,P_1,\ldots,P_t=Q$ such that $P_{i-1}\cap P_i$ is a ridge for every $i=1,\ldots,t$ (see \cite[11.6]{Bjo95}). It is \emph{normal} if $\operatorname{lk}_\Delta F$ is connected whenever $\dim\operatorname{lk}_\Delta F\geq1$. A \emph{pseudomanifold} is a pure strongly connected complex in which every ridge belongs to exactly two facets. Thus a pseudomanifold is normal if and only if its Stanley--Reisner ring satisfies $(S_2)$.

There is a similar description of $(G_2)$. The ring $\kk[\Delta]$ satisfies $(G_2)$ if and only if, for every $F\in\Delta$, $\operatorname{lk}_\Delta F$ consists of one or two vertices when $\dim\operatorname{lk}_\Delta F=0$, and is a cycle or a path with one or two edges when $\dim\operatorname{lk}_\Delta F=1$. Indeed, for $F\in\Delta$, set $\mathfrak p_F=(x_i:i\notin F)R$ and $K=\kk(x_i:i\in F)$. If $\mathfrak n$ denotes the graded maximal ideal of $K[\operatorname{lk}_\Delta F]$, then $R_{\mathfrak p_F}\cong K[\operatorname{lk}_\Delta F]_{\mathfrak n}$. The assertion follows from these localizations and the Gorenstein criterion for Stanley--Reisner rings (see \cite[Theorem~5.6.1]{BH98}).

\begin{lemma}\label{lem:g2-pseudomanifold}
Let $\Delta=\core(\Delta)$ have dimension at least one, and suppose that $R=\kk[\Delta]$ satisfies $(S_2)$. Then $R$ satisfies $(G_2)$ if and only if $\Delta$ is a normal pseudomanifold.
\end{lemma}
\begin{proof}
The condition $(S_2)$ implies that $\Delta$ is pure and normal, and hence strongly connected by \cite[proof of Proposition~11.7]{Bjo95}. In dimension at least two, $(G_2)$ implies that $\Delta$ is a pseudomanifold by \cite[Proposition~4.13]{KM26}. In dimension one, the link criterion gives a cycle or a path with one or two edges; the paths are excluded because they have cone vertices. Conversely, if $\Delta$ is a normal pseudomanifold, then $\operatorname{lk}_\Delta F$ consists of two vertices when $\dim\operatorname{lk}_\Delta F=0$, and is a cycle when $\dim\operatorname{lk}_\Delta F=1$. The link criterion therefore gives $(G_2)$.
\end{proof}

We write $\widetilde H_i(-;\kk)$ for reduced simplicial homology, computed from augmented chains. In particular, $\widetilde H_{-1}(\{\varnothing\};\kk)=\kk$. A pure $(d-1)$-dimensional complex is a \emph{$\kk$-homology sphere} if, for every face $F$, including the empty face, $\widetilde H_i(\operatorname{lk}_\Delta F;\kk)$ is $\kk$ for $i=d-|F|-1$ and zero otherwise. A \emph{$\kk$-homology manifold} is a connected pure complex whose nonempty face links are $\kk$-homology spheres.

For a squarefree module $M$, we write $M_F$ for $M_{\chi_F}$. If $\mathbf{u}\in\NN^n$ has support $F$, multiplication by $\mathbf{x}^{\mathbf{u}-\chi_F}$ is an isomorphism from $M_F$ to $M_{\mathbf{u}}$. Thus the squarefree components determine all homogeneous components of $M$. In particular, $M$ is generated in squarefree multidegrees. The canonical module of a squarefree module is again squarefree (see \cite[Section~2 and Theorem~2.6]{Yan00}).

We shall use the following description due to Gr\"abe \cite{Gra84}. Let $R=\kk[\Delta]$ have dimension $d$, and fix orientations of the facets of dimension $d-1$. For $F\in\Delta$, the component $[\omega_R]_F$ consists of the linear combinations of these facets containing $F$ whose boundaries have coefficient zero at every $(d-2)$-face containing $F$. In particular,
\[
\dim_\kk[\omega_R]_F=\dim_\kk\widetilde H_{d-|F|-1}(\operatorname{lk}_\Delta F;\kk).
\]

This description also gives the multiplication maps. If $F\subseteq G$ are faces of $\Delta$, multiplication by $x_{G\setminus F}$ retains the coefficients of the facets containing $G$ and discards the others.

For a pure complex $\Delta$ and a face $F$, set $A=\kk[\operatorname{lk}_\Delta F]$. Under the usual identification $R_{x_F}\cong A[x_i^{\pm1}:i\in F]$, flat base change for $\Ext$ identifies $(\omega_R)_{x_F}$ with $\omega_A\otimes_A R_{x_F}$ as ungraded modules; purity ensures that $\dim A=\dim R-|F|$. Thus \cite[Lemma~1.5\,(iii)]{HHS} gives
\begin{equation}\label{eq:trace-link}
\tr_R(\omega_R)R_{x_F}=\tr_A(\omega_A)R_{x_F}.
\end{equation}

\section{Proofs of the main results}\label{sec:trace-structure}

In this section we prove Theorems~\ref{thm:main-a} and~\ref{thm:main-b}. Throughout, $S=\kk[x_1,\ldots,x_n]$ denotes the polynomial ring with $\deg x_i=\ee_i$. We first bound the degrees of generators of trace ideals, with applications to canonical modules and exterior powers of the module of differentials. We then determine when the canonical trace ideal is the square of its radical.

Miyazaki showed that a module generated in a single multidegree over a Stanley--Reisner ring has a radical trace ideal in his Japanese notes for the 37th Seminar on Commutative Ring Theory (see \cite[p.~8]{MiyazakiSem26}). His argument lowers the degree of a homomorphism by dividing its values by a variable, and the same argument also gives the cube-free case. The following lemma extends this argument to homomorphisms with values in an arbitrary positively determined module.

For $S$-modules $M$ and $N$, the \emph{trace module of $M$ in $N$} is
\[
\tau_M(N)=\sum_{\varphi\in\Hom_S(M,N)}\varphi(M)
\]
(see \cite[Definition~2.1]{LP22}). In particular, $\tau_M(S)=\tr_S(M)$.

\begin{lemma}\label{lem:images-bound}
Let $\mathbf{a},\mathbf{b}\in\NN^n$ with $\mathbf{a}\geq\one$, and let $N$ be a positively $\mathbf{a}$-determined $S$-module. Suppose that, after a shift, the multigraded $S$-module $M$ is generated in degrees belonging to $[0,\mathbf{b}]$. Then the trace module $\tau_M(N)$ is positively $(\mathbf{a}+\mathbf{b})$-determined.
\end{lemma}
\begin{proof}
Set $T=\tau_M(N)$. After shifting $M$, choose homogeneous generators in degrees belonging to $[0,\mathbf{b}]$. Fix $i$ and $\mathbf{u}\in\NN^n$ with $u_i\geq a_i+b_i$. The space $T_{\mathbf{u}+\ee_i}$ is spanned by values $\varphi(m)$, where $m$ is one of these generators and $\varphi$ is homogeneous of degree $\mathbf{u}+\ee_i-\deg m$. Thus $(\deg\varphi)_i\geq a_i+1$. Since $M$ is nonnegatively graded, $\varphi(M)$ lies in $\bigoplus_{v_i\geq a_i+1}N_{\mathbf{v}}$. As $N$ is positively $\mathbf{a}$-determined, multiplication by $x_i$ induces an $S$-linear isomorphism
\[
\bigoplus_{v_i\geq a_i}N_{\mathbf{v}}\xrightarrow{\ x_i\ }\bigoplus_{v_i\geq a_i+1}N_{\mathbf{v}}.
\]
Composing $\varphi$ with its inverse gives a homogeneous map $\psi$ with $\varphi=x_i\psi$. Hence $T_{\mathbf{u}+\ee_i}=x_iT_{\mathbf{u}}$. Multiplication by $x_i$ is injective on $T_{\mathbf{u}}$ because $T\subseteq N$ and $u_i\geq a_i$. This proves the assertion.
\end{proof}

The following theorem is our first main result.

\begin{theorem}\label{thm:cube-free}
Let $I\subsetneq S$ be a monomial ideal, and set $R=S/I$. Let $\mathbf{a},\mathbf{b}\in\NN^n$ with $\mathbf{a}\geq\one$, and suppose that $R$ is positively $\mathbf{a}$-determined. If a multigraded $R$-module $M$ is generated, after a shift, in degrees belonging to $[0,\mathbf{b}]$, then $\tr_R(M)$ is $(\mathbf{a}+\mathbf{b})$-bounded. In particular, when $\mathbf{b}\geq\one$, the conclusion holds if a shift of $M$ is positively $\mathbf{b}$-determined.
\end{theorem}
\begin{proof}
Apply Lemma~\ref{lem:images-bound} with $N=R$, noting that $\Hom_R(M,R)=\Hom_S(M,R)$. It follows that $\tr_R(M)$ is positively $(\mathbf{a}+\mathbf{b})$-determined, and hence generated by monomials whose degrees belong to $[0,\mathbf{a}+\mathbf{b}]$.
\end{proof}

\begin{corollary}\label{cor:monomial-canonical-bound}
Let $I\subsetneq S$ be an $\mathbf{a}$-bounded monomial ideal, where $\mathbf{a}\geq\one$, and set $R=S/I$. If $L$ is a positively $\mathbf{b}$-determined $R$-module with $\mathbf{b}\geq\one$, then $\tr_R(\Ext_S^q(L,\omega_S))$ is $(\mathbf{a}+\mathbf{b})$-bounded for every $q\geq0$. In particular, $\tr_R(\omega_R)$ is $2\mathbf{a}$-bounded.
\end{corollary}
\begin{proof}
Set $E=\Ext_S^q(L,S(-\mathbf{b}))$. By Remark~\ref{rem:positively-determined}, $E$ is positively $\mathbf{b}$-determined. Moreover, $I$ annihilates $E$, so $E$ is an $R$-module. Since $E$ differs from $\Ext_S^q(L,\omega_S)$ only by a shift, the asserted bound follows from Theorem~\ref{thm:cube-free}. For the last assertion, take $L=R$, $\mathbf{b}=\mathbf{a}$, and $q=n-\dim R$.
\end{proof}

\begin{proof}[Proof of Theorem~\ref{thm:main-a}]
By Remark~\ref{rem:positively-determined}, $R$ is positively $\mathbf{a}$-determined and $M$ is generated in multidegrees belonging to $[0,\mathbf{b}]$. The first assertion therefore follows from Theorem~\ref{thm:cube-free}. The second assertion is the last statement of Corollary~\ref{cor:monomial-canonical-bound}.
\end{proof}

\begin{remark}
Taking $\mathbf{a}=\one$ and $\mathbf{b}=0$ in Theorem~\ref{thm:cube-free}, we recover the fact that the trace ideal of a module generated in a single multidegree over a Stanley--Reisner ring is radical. If $\Delta$ is doubly Cohen--Macaulay, then $\omega_R$ is generated in multidegree zero (see \cite[Corollary~2.6]{Yan11}). Thus $\kk[\Delta]$ has radical canonical trace ideal, as proved by Miyazaki \cite[Theorem~3.4]{MiyazakiSem26}.

Following Miyazaki \cite[Definition~3.2]{Miyazaki}, a Cohen--Macaulay ring $A$ has the \emph{canonical trace radical (CTR) property} if, for every $\mathfrak p\in\operatorname{Spec}A$, the local ring $A_{\mathfrak p}$ has a canonical module whose trace ideal is radical. For a Cohen--Macaulay $\NN$-graded $\kk$-algebra $A$ with $A_0=\kk$, this is equivalent to $\tr_A(\omega_A)$ being radical (see \cite[Proposition~3.4]{Miyazaki}).
\end{remark}

Theorem~\ref{thm:cube-free} also applies to exterior powers of the module of differentials, whose trace ideals were studied in \cite{IM26}.

\begin{corollary}\label{cor:differentials-cube-free}
Let $I\subsetneq S$ be an $\mathbf{a}$-bounded monomial ideal with $\mathbf{a}\geq\one$, and set $R=S/I$. For every $p\geq0$, the trace ideal $\tr_R(\bigwedge_R^p\Omega_{R/\kk})$ is $(\mathbf{a}+\one)$-bounded. In particular, it is cube-free when $R$ is a Stanley--Reisner ring.
\end{corollary}
\begin{proof}
Give $\Omega_{R/\kk}$ the multigrading $\deg dx_i=\ee_i$. In the usual presentation (see \cite[p.~195, Example]{Mat86}), the generators are $dx_1,\ldots,dx_n$, and the relations are $d(\mathbf{x}^{\mathbf{u}})=0$ for the minimal monomial generators $\mathbf{x}^{\mathbf{u}}$ of $I$. Each relation $d(\mathbf{x}^{\mathbf{u}})=\sum_{i:u_i>0}u_i\mathbf{x}^{\mathbf{u}-\ee_i}\,dx_i$ is homogeneous of degree $\mathbf{u}$. The exterior products $dx_{i_1}\wedge\cdots\wedge dx_{i_p}$ therefore generate $\bigwedge_R^p\Omega_{R/\kk}$ in degrees $\chi_F$, where $F=\{i_1,\ldots,i_p\}$. The conclusion follows from Theorem~\ref{thm:cube-free} with $\mathbf{b}=\one$.
\end{proof}

\begin{remark}
The bound $2\mathbf{a}$ in Corollary~\ref{cor:monomial-canonical-bound} is sharp in every coordinate. Let $\Delta$ be a six-vertex triangulation of $\mathbb{RP}^2$, assume that $\operatorname{char}\kk\ne2$, and set $S=\kk[x_1,\ldots,x_6]$ and $A=S/I_\Delta$. By \cite[Corollary~3.10]{MV}, we have $\tr_A(\omega_A)=\mathfrak m_A^2$.

Fix $\mathbf{a}\in\NN^6$ with $\mathbf{a}\geq\one$, and set $T=\kk[y_1,\ldots,y_6]$. Consider the $\kk$-algebra homomorphism $S\to T$ given by $x_i\mapsto y_i^{a_i}$. Each monomial in $T$ has a unique expression $(y_1^{a_1})^{v_1}\cdots(y_6^{a_6})^{v_6}y_1^{u_1}\cdots y_6^{u_6}$ with $0\leq u_i<a_i$. Hence $T$ is a free $S$-module with basis $\{y_1^{u_1}\cdots y_6^{u_6}:0\leq u_i<a_i\}$. Set $B=T/I_\Delta T$. Thus $B\cong A\otimes_S T$ is a finite free $A$-module, and $\dim B=\dim A=3$. The ideal $I_\Delta T$ is obtained by replacing $x_i$ by $y_i^{a_i}$ in the squarefree generators of $I_\Delta$, so it is $\mathbf{a}$-bounded.

For this argument, give $A$ the grading $\deg x_i=a_i$, so that $A\to B$ is graded. The ring $A$ is Cohen--Macaulay, and $B/\mathfrak m_AB\cong\kk[y_1,\ldots,y_6]/(y_1^{a_1},\ldots,y_6^{a_6})$ is an Artinian complete intersection. Applying \cite[Theorem~3.3.14\,(a)]{BH98} at the homogeneous maximal ideals and using graded uniqueness (see \cite[Proposition~3.6.9]{BH98}), we obtain $\omega_B\cong\omega_A\otimes_A B$ as ungraded $B$-modules.

Since $B$ is flat over $A$ and $\mathfrak m_AB=(y_1^{a_1},\ldots,y_6^{a_6})B$, flat base change for trace ideals (see \cite[Lemma~1.5\,(iii)]{HHS}) now gives
\[
\tr_B(\omega_B)=\tr_A(\omega_A)B=\mathfrak m_A^2B=(y_1^{a_1},\ldots,y_6^{a_6})^2B.
\]
The defining ideal of $B$ contains no pure power of a variable. It follows that each $y_i^{2a_i}$ is a minimal generator of this trace ideal.
\end{remark}

For normal pseudomanifolds, Gr\"abe's description takes the following form.

\begin{lemma}\label{lem:link-cycles}
Let $\Delta$ be a normal pseudomanifold of positive dimension, and set $R=\kk[\Delta]$. For every face $F$ with $[\omega_R]_F\ne0$, the space $[\omega_R]_F$ is one-dimensional, with a basis element $g_F$ whose coefficients at all facets containing $F$ are $\pm1$. If $F\subseteq H\in\Delta$, then $x_{H\setminus F}g_F=\pm g_H$.
\end{lemma}
\begin{proof}
By \cite[Lemma~3.2]{MV}, each nonzero component is one-dimensional, and multiplication by $x_{H\setminus F}$ gives an isomorphism to the component at $H$. The proof of \cite[Lemma~2.5]{MV} allows us to choose each $g_F$ with coefficients $\pm1$. Multiplication retains the coefficients at facets containing $H$, so $x_{H\setminus F}g_F=\pm g_H$.
\end{proof}

We extend the pairing for homology manifolds in \cite[proof of Proposition~3.9]{MV} to normal pseudomanifolds and show that it identifies the canonical module with its $R$-dual.

\begin{lemma}\label{lem:canonical-pairing}
Let $\Delta$ be a normal pseudomanifold of positive dimension, and set $R=\kk[\Delta]$ and $\mathfrak a=(x_F:F\in\Delta,\ [\omega_R]_F\ne0)R$. Then $\tr_R(\omega_R)=\mathfrak a^2$ and $\sqrt{\tr_R(\omega_R)}=\mathfrak a$.
\end{lemma}
\begin{proof}
Set $M=\omega_R$ and $d=\dim R$. Fix orientations of the facets, and set $B=\prod_{P\text{ facet}}R/\mathfrak p_P$. The natural maps embed $R$ in $B$. Gr\"abe's description of multiplication and squarefreeness give a multigraded $R$-linear embedding $M\hookrightarrow B$ determined on $M_F$ by $\sum_{P\supseteq F}c_PP\mapsto(c_Px_F)_P$, where $c_P=0$ for $P\not\supseteq F$. We identify $M$ with its image.

Choose the bases $g_F$ from Lemma~\ref{lem:link-cycles}. If $M_F,M_G\ne0$ and $H=F\cup G$ is a face, then $x_{H\setminus F}g_F$ and $x_{H\setminus G}g_G$ are $\pm g_H$. The product of their coefficients is therefore constant on the facets containing $H$, since the coefficients of $g_H$ have square one. At all other facets the products vanish. Hence multiplication in $B$ gives $g_Fg_G=\pm x_Fx_G$. If $F\cup G$ is not a face, the product is zero. Since the $g_F$ generate $M$, multiplication in $B$ defines a symmetric $R$-bilinear pairing $b:M\times M\to R$ of multidegree zero. The induced map $M\otimes_R M\to R$ has image $\mathfrak a^2$.

We claim that the multigraded map $\Phi:M\to\Hom_R(M,R)$ given by $t\mapsto b(t,-)$ is an isomorphism. The canonical module $M$ satisfies $(S_2)$ (see \cite[Lemma~1.9\,(e)]{Sch98}), and $\Hom_R(M,R)$ satisfies $(S_1)$ because it embeds in a finite free module over the reduced ring $R$. By \cite[Lemma~2.20]{DFI15}, it suffices to check $\Phi$ at primes of height at most one.

Let $\mathfrak p$ be such a prime, and set $F=\{i:x_i\notin\mathfrak p\}$. Then $F\in\Delta$ and $\mathfrak p_F\subseteq\mathfrak p$. Purity gives $\operatorname{ht}\mathfrak p_F=d-|F|\leq1$, so $F$ is a facet or a ridge. After localizing at $x_F$, Gr\"abe's description gives $M_{x_F}=R_{x_F}g_F$, and $b(g_F,g_F)=x_F^2$ is a unit in $R_{x_F}$. Thus $\Phi_{x_F}$, and hence $\Phi_{\mathfrak p}$, is an isomorphism. The surjectivity of $\Phi$ shows that $\tr_R(M)$ is generated by the values of $b$, and hence equals $\mathfrak a^2$. Since $\mathfrak a$ is generated by squarefree monomials, it is radical, and the last assertion follows.
\end{proof}

The next lemma proves the converse in Theorem~\ref{thm:main-b} assuming only that $\Delta$ is pure.

\begin{lemma}\label{lem:low-dimensional-trace}
Let $\Delta$ be pure, and set $R=\kk[\Delta]$. If $R$ does not satisfy $(G_2)$, there is a face $F$ such that $\mathfrak p_F$ has height one or two and $\tr_R(\omega_R)R_{\mathfrak p_F}=\mathfrak p_FR_{\mathfrak p_F}\ne0$. Consequently, if the canonical trace ideal of $R$ is the square of its radical, then $R$ satisfies $(G_2)$.
\end{lemma}
\begin{proof}
Choose $F$ maximal among the faces for which $\Gamma=\operatorname{lk}_\Delta F$ has dimension at most one and $\kk[\Gamma]$ is not Gorenstein. Such a face exists by the link criterion in Section~\ref{sec:preliminaries}. Set $A=\kk[\Gamma]$, and let $\mathfrak n$ denote its graded maximal ideal.

If $\dim\Gamma=1$, purity and the choice of $F$ imply that every vertex has degree one or two, so each connected component is a path or a cycle. If $\Gamma$ is connected, it is a path with at least three edges, and \cite[Theorem~4.19\,(3)]{KM26} gives $\tr_A(\omega_A)=\mathfrak n$. If $\Gamma$ is disconnected, the same equality follows from \cite[Theorem~4.21\,(2)]{KM26}. This also covers the case $\dim\Gamma=0$, in which $\Gamma$ has at least three vertices.

Set $T=\tr_R(\omega_R)$. By \eqref{eq:trace-link}, we have $TR_{x_F}=\mathfrak nR_{x_F}=\mathfrak p_FR_{x_F}$. Localizing further at $\mathfrak p_F$ gives $T_{\mathfrak p_F}=\mathfrak p_FR_{\mathfrak p_F}$. Since $\dim R_{\mathfrak p_F}=\operatorname{ht}\mathfrak p_F=\dim A$ is one or two, $R_{\mathfrak p_F}$ is not a field. Its maximal ideal $\mathfrak p_FR_{\mathfrak p_F}$ is therefore nonzero. If $T=(\sqrt T)^2$, localization would give $\mathfrak p_FR_{\mathfrak p_F}=(\mathfrak p_FR_{\mathfrak p_F})^2$, contrary to Nakayama's lemma.
\end{proof}

We now prove the second main result of this paper.

\begin{proof}[Proof of Theorem~\ref{thm:main-b}]
The assertion is immediate for $R=0$. Suppose that $R\ne0$ satisfies $(S_2)$ and $(G_2)$. Set $A=\kk[\core(\Delta)]$, so that $R$ is a polynomial ring over $A$. By flat base change for $\Ext$ and trace ideals (see \cite[Lemma~1.5\,(iii)]{HHS}), we have $\tr_R(\omega_R)=\tr_A(\omega_A)R$. Radicals also commute with polynomial extension. Since $(S_2)$ and $(G_2)$ are preserved on passing between $A$ and $R$, we may assume that $\Delta=\core(\Delta)$.

If $\dim\Delta\leq0$, then $(G_2)$ implies that $R$ is Gorenstein, and the equality follows. Otherwise Lemma~\ref{lem:g2-pseudomanifold} shows that $\Delta$ is a normal pseudomanifold, and Lemma~\ref{lem:canonical-pairing} gives the equality.

Conversely, $(S_2)$ implies that $\Delta$ is pure, so Lemma~\ref{lem:low-dimensional-trace} applies.
\end{proof}

\begin{corollary}\label{cor:core-square}
Let $\Delta=\core(\Delta)$ be a simplicial complex of dimension at least one, and let $R=\kk[\Delta]$ satisfy $(S_2)$. Then $\Delta$ is a normal pseudomanifold if and only if
\[
\tr_R(\omega_R)=\bigl(\sqrt{\tr_R(\omega_R)}\bigr)^2.
\]
\end{corollary}
\begin{proof}
Apply Theorem~\ref{thm:main-b} and Lemma~\ref{lem:g2-pseudomanifold}.
\end{proof}

\begin{corollary}\label{cor:radical-qg}
Let $R=\kk[\Delta]\ne0$ satisfy $(S_2)$ and $(G_2)$. Then $\tr_R(\omega_R)$ is radical if and only if $R$ is quasi-Gorenstein.
\end{corollary}
\begin{proof}
By \cite[Remark~2.13]{KM26}, $\omega_R$ is isomorphic to an ideal containing a nonzerodivisor, so $T=\tr_R(\omega_R)$ also contains a nonzerodivisor. If $T$ is radical, Theorem~\ref{thm:main-b} gives $T=T^2$, whence $T=R$ by the graded Nakayama lemma. Now apply \cite[Remark~2.9\,(4)]{KM26}.
\end{proof}

\begin{example}\label{ex:square-nonisolated}
Let $\Delta$ be a six-vertex triangulation of $\mathbb{RP}^2$, and assume that $\operatorname{char}\kk\ne2$. Let $\Gamma$ consist of two isolated vertices $p,q\notin V(\Delta)$, and set $\Sigma=\Delta*\Gamma$. Thus $\Sigma$ is the suspension of $\Delta$. This is the construction used by Miyashita and Varbaro (see \cite[Remark~2.7]{MV}), with the torus replaced by $\Delta$.

Set $R=\kk[\Sigma]$ and $\mathfrak a=(x_i:i\in V(\Delta))R$. Then $R$ is Cohen--Macaulay and $\Sigma$ is a normal pseudomanifold. Since $\operatorname{lk}_\Sigma\{p\}=\operatorname{lk}_\Sigma\{q\}=\Delta$, $\Sigma$ is not a $\kk$-homology manifold. We have $\widetilde H_2(\Delta;\kk)=0$, whereas $\operatorname{lk}_\Delta F$ is a $\kk$-homology sphere for every nonempty face $F\in\Delta$. Gr\"abe's description and Lemma~\ref{lem:canonical-pairing} therefore give $\tr_R(\omega_R)=\mathfrak a^2$.
\end{example}

\begin{example}\label{ex:square-hypotheses}
The following two examples show that $(S_2)$ cannot simply be omitted from Theorem~\ref{thm:main-b} and that $(G_2)$ cannot be replaced by $(G_1)$. Both examples work over an arbitrary field.

First, let $A=\kk[x_1,x_2,x_3]$, $B=\kk[y_1,y_2,y_3]$, and $D=A\times_\kk B$. Set $R=D[p,q]/(pq)$ and $\mathfrak a=(x_1,x_2,x_3,y_1,y_2,y_3)R$. Let $\Delta$ denote the defining complex of $R$, which is the suspension of two disjoint $2$-simplices. It is connected and pure, and it has no cone vertices. By \cite[Theorem~3.11]{KM26}, $\tr_D(\omega_D)=\mathfrak m_D$. Since $\kk[p,q]/(pq)$ is isomorphic to its canonical module, we have $\omega_R\cong\omega_D\otimes_\kk\kk[p,q]/(pq)$. Flat base change for trace ideals now gives $\tr_R(\omega_R)=\mathfrak a$. This ideal is radical and differs from its square, so the equality in Theorem~\ref{thm:main-b} fails.

For every face $F$ of this complex, $\operatorname{lk}_\Delta F$ consists of one or two vertices when $\dim\operatorname{lk}_\Delta F=0$, and is a path with one or two edges when $\dim\operatorname{lk}_\Delta F=1$. Thus $R$ satisfies $(G_2)$ by the link criterion, whereas $\operatorname{lk}_\Delta\{p\}$ is disconnected and has dimension two, so $R$ does not satisfy $(S_2)$.

Next, let $\Delta$ be the path with facets $\{1,2\}$, $\{2,3\}$, and $\{3,4\}$, and set $R=\kk[\Delta]$. Then $R$ is Cohen--Macaulay of dimension two, and \cite[Theorem~A]{MV} gives $\tr_R(\omega_R)=\mm\ne\mm^2$. Its non-Gorenstein locus is $\{\mm\}$, so $R$ satisfies $(S_2)$ and $(G_1)$ but not $(G_2)$.
\end{example}

\medskip
\noindent\textbf{Acknowledgments.} The author thanks Naoki Terai for helpful discussions on Stanley--Reisner rings and for his encouragement. The author was supported by JSPS KAKENHI Grant Number JP25KJ1744. OpenAI's ChatGPT was used for preliminary mathematical discussions, literature searches, language editing, and \TeX{} preparation. The author takes responsibility for the mathematical content, references, and final text.
\par

\end{document}